\documentclass[12pt,reqno]{amsart}
\usepackage{amsfonts}
\usepackage{amsmath,amsthm,amsfonts,amssymb}
\newtheorem{lemma}{Lemma}[section]
\newtheorem{theorem}{Theorem}[section]
\newtheorem{corollary}{Corollary}[section]
\newtheorem*{Main Theorem}{Main Theorem}

\def\bl{\begin{lemma}}
\def\bt{\begin{theorem}}
\def\bbt{\begin{Main Theorem}}
\def\el{\end{lemma}}
\def\et{\end{theorem}}
\def\eet{\end{Main Theorem}}
\def\bp{\begin{proof}}
\def\ep{\end{proof}}
\def\bc{\begin{corollary}}
\def\ec{\end{corollary}}

\def\mb{\mathbb}

\def\O{\Omega}
\def\a{\alpha}

\def\p{\partial}

\def\-{\setminus}

\def\vp{\varphi}

\def\lt{\left}
\def\rt{\right}
\def\+{\bigcup}
\def\.{\bigcap}

\title[]
{$L^2$ estimates of Poincar\'e-Lelong equations on convex domains in $\mathbb{C}^n$}

\author []{Shaoyu Dai$^1$, Yang Liu$^2$ and Yifei Pan$^3$}

\address{1 Department of Mathematics, Jinling Institute of Technology, Nanjing, 211169, China.}
\address{\it E-mail address: dymdsy@163.com}

\address{2 Department of Mathematics, Zhejiang Normal University, Jinhua, 321004, China.}
\address{\it E-mail address: liuyang@zjnu.edu.cn}

\address{3 Department of Mathematical Sciences, Purdue University Fort Wayne, Fort Wayne, 46805-1499, USA.}
\address{\it E-mail address: pan@pfw.edu}

\begin{document}

\begin{abstract}
In this paper, we prove the existence of solutions of the Poincar\'e-Lelong equation $\sqrt{-1}\partial\overline{\partial}u=f$
on a strictly convex bounded domain $\Omega\subset\mathbb{C}^n$
$(n\geq1)$, where $f$ is a $d$-closed $(1,1)$ form and is in the weighted Hilbert space $L^2_{(1,1)}(\Omega,e^{-\varphi})$. The novelty of this paper is to apply a weighted $L^2$ version of Poincar\'e Lemma for real $2$-forms, and then apply H\"{o}rmander's $L^2$ solutions for Cauchy-Riemann equations.
\end{abstract}

\maketitle

\section{Introduction}

In this paper, a continuation of \cite {we}, we will study the Poincar\'e-Lelong equation $\sqrt{-1}\partial\overline{\partial}u=f$ in a weighted Hilbert space on a strictly convex bounded domain $\O\subset\mb{C}^n$ $(n\geq1)$. Using a weighted $L^2$ version of Poincar\'e Lemma for real forms, we obtain the existence of solutions with the norm estimate. More precisely, we prove the following theorem.

\bbt
Let $\O$ be a strictly convex bounded domain in $\mb{C}^n$. Let $\vp$ be a nonnegative smooth function on $\overline{\O}$ such that
$$\sum^{2n}_{j,k=1}
\frac{\p^2\vp}{\p x_j\p x_k}\xi_{j}\xi_{k}\geq c|\xi|^2\ \ \mbox{for all}\ \  \xi=(\xi_1,\cdots,\xi_{2n})\in\mathbb{R}^{2n},$$
where $c>0$ is a constant.
Then, for each $(1,1)$ form $f$ in the weighted Hilbert space
$L^2_{(1,1)}(\O,e^{-\vp})$ with $\p f=\overline{\partial}f=0$, there exists a solution $u$ in $L^2(\O,e^{-\vp})$ solving the Poincar\'e-Lelong equation
$$\sqrt{-1}\partial\overline{\partial}u=f$$ in $\O$, in the sense of distributions,
with the norm estimate
\begin{align*}
\int_{\O} |u|^2e^{-\vp}\leq
\frac{8}{c^2}\int_{\O}|f|^2e^{-\vp}.
\end{align*}
\eet

\bigskip
\noindent\textbf{Corollary.}
\textit{
For each $(1,1)$ form $f$ in $L^2_{(1,1)}(\O)$ with $\p f=\overline{\partial}f=0$, there exists a solution $u$ in $L^2(\O)$ solving the Poincar\'e-Lelong equation
$\sqrt{-1}\partial\overline{\partial}u=f$ in $\O$, in the sense of distributions,
with the $L^2$ norm estimate
\begin{align*}
\int_{\O} |u|^2\leq
c_\O\int_{\O}|f|^2,
\end{align*}
where $c_\O>0$ is a constant depended on the diameter of $\O$.
}

\bigskip
For the Poincar\'e-Lelong equation $\sqrt{-1}\partial\overline{\partial}u=f$, P. Lelong \cite{1} studied it in connection with questions on entire functions, and showed, unexpectedly, that with suitable restrictions on the growth of $f$,
the equation could be reduced to solving the more familiar equation $\frac{1}{4}\Delta u=\mbox{trace}(f)$ (Poisson equation). Mok, Siu and Yau \cite{2} studied the equation on a complete K\"{a}hler manifold and obtained important applications to questions on when a (noncompact) K\"{a}hler manifold is biholomorphicly equivalent to $\mb{C}^n$.
Andersson \cite{Mats Andersson} studied the Poincar\'e-Lelong equation for smooth forms in the unit ball using integral representations. Recently, Chen \cite{3} obtained solutions of the equation when $f$ is assumed to be a smooth $(1,1)$ $d$-closed form with compact support in $\mb{C}^n$, and he applied his result to prove a version of Hartog's extension theorem for pluriharmonic functions (for related results see also \cite{L}, \cite{dbar} and \cite{ddd}).

Recently, we \cite{we} studied the Poincar\'e-Lelong equation in the whole space $\mb{C}^n$, and proved the existence of (global) solutions in the weighted Hilbert space with Gaussian measure as follows.

\bigskip
\noindent\textbf{Theorem.}
\textit{
For each $(1,1)$ form $f$ in the weighted Hilbert space $ L^2_{(1,1)}(\mb{C}^n,e^{-|z|^2})$ with $\p f=\overline{\partial}f=0$, there exists a solution $u$ in $L^2(\mb{C}^n,e^{-|z|^2})$ solving the Poincar\'e-Lelong equation
$$\sqrt{-1}\partial\overline{\partial}u=f$$ in $\mb{C}^n$, in the sense of distributions,
with the norm estimate $$\int_{\mb{C}^n} |u|^2e^{-|z|^2}\leq
2\int_{\mb{C}^n} |f|^2e^{-|z|^2}.$$
}

As a matter of fact, the key idea of the proof of the main theorem is quite similar to that of the theorem above. First we convert the $(1,1)$ form $f$ to a real $2$-form. Second for the real $2$-form, we apply a weighted $L^2$ version of Poincar\'e Lemma that we shall give a detailed proof. At last, we apply H\"{o}rmander's $L^2$ solutions for the Cauchy-Riemann equations. Since the domain considered has a smooth boundary, we have to study carefully about the adjoint of the $\p\overline{\p}$ operator following the approach of Berndtsson \cite{B}, who essentially gave the second proof of the H\"{o}rmander's theorem for $\overline{\p}$.

It was Berndtsson \cite{5}, who first studied the $d$-equation for real $1$-forms with morse function weights and pointed out that
the H\"{o}rmander's $L^2$ method could be used for the $d$-equation in convex domains and with a convex weight function. Since our proof of the main theorem depends significantly on a weighted $L^2$ version of Poincar\'e Lemma (below), and the classical Poincar\'e Lemma would not provide a $L^2$ estimates for $d$-equation, we decide to include a detailed proof of Poincar\'e Lemma. In addition, the proof will provide a specific constant that we shall use in the main theorem.

\bigskip
\noindent\textbf{Poincar\'e Lemma.} (A weighted $L^2$ version for $p+1$-forms)
\textit{
Let $N\geq1$ be an integer and $G$ be a strictly convex bounded domain in $\mathbb{R}^N$. Let $\vp$ be a nonnegative smooth function on $\overline{G}$ such that
$$\sum^N_{j,k=1}
\frac{\p^2\vp}{\p x_j\p x_k}\xi_{j}\xi_{k}\geq c|\xi|^2\ \ \mbox{for all}\ \ \xi=(\xi_1,\cdots,\xi_{N})\in\mathbb{R}^{N},$$
where $c$ is a positive constant.
Let $p$ be an integer with $0\leq p\leq N-1$.
Then, for each $f$, a $d$-closed $p+1$-form in the weighted Hilbert space $ L^2_{p+1}(G,e^{-\vp})$, there exists a solution $p$-form $u$ in $L^2_{p}(G,e^{-\vp})$ solving equation
$$du=f$$ in $G$, in the sense of distributions,
with the norm estimate $$\int_G |u|^2e^{-\vp}\leq
\frac{1}{c(p+1)}\int_G|f|^2e^{-\vp}.$$
}

In this paper, Section 2 and 3 are for Poincar\'e Lemma; Section 4 and 5 are for the main theorem.

\section{Preliminary for Poincar\'e Lemma}

Let $N\geq1$ and $p\geq0$ be integers.
For multiindex $I=(i_1,\cdots,i_p)$, where $i_1,\cdots,i_p$ are integers between $1$ and $N$, define $|I|=p$ and $dx^I=dx_{i_1}\wedge\cdots\wedge dx_{i_p}$.
Let $G$ be a strictly convex bounded domain in $\mb{R}^N$.
In general, a $p$-form $f$ on $G$ is a formal combination
$$f=\sum'_{|I|=p}f_Idx^I,$$
where ${\sum}'$ implies that the summation is performed only over strictly increasing multi-indices and
$f_I:G \rightarrow\mb{R}$ is a real-valued function for each $I$.
For $p$-forms $f$ and $g$, we denote their pointwise scalar product by
$$f\cdot g=\sum'_{|I|=p}f_Ig_I.$$
Let $\vp$ be a nonnegative smooth function on $\overline{G}$ and the weighted Hilbert space for $p$-forms
$$L^2_p(G,e^{-\vp})
=\{f=\sum'_{|I|=p}f_Idx^I\mid f_I\in L^2_{loc}(G); \int_{G}|f|^2 e^{-\vp}<+\infty\},$$
where $|f|^2=f\cdot f$. We denote
the weighted inner product for $f,g\in L^2_p(G,e^{-\vp})$ by
$$\langle f,g\rangle_{L^2_p(G,e^{-\vp})}=\int_{G}
f\cdot g e^{-\vp},$$
and the weighted norm of $f\in L^2_p(G,e^{-\vp})$ by
$\|f\|_{L^2_p(G,e^{-\vp})}=\sqrt{\langle f,f\rangle}_{L^2_p(G,e^{-\vp})}.$
In particular, we denote
$$L^2_0(G,e^{-\vp})=L^2(G,e^{-\vp})
=\{f\mid f\in L^2_{loc}(G); \int_{G}|f|^2 e^{-\vp}<+\infty\}.$$

For simplicity, we will write $L^2_p(e^{-\vp})$ for $L^2_p(G,e^{-\vp})$ in Section 2 and 3, since we only deal with $G$ in these sections.

Let $\mathcal{D}_p$ denote the set of $p$-forms whose coefficients are smooth functions with compact support in $G$.
For $p$-form $u={\sum'_{|I|=p}}u_Idx^I\in L^2_p(e^{-\vp})$, in the sense of distributions, the differential of $u$ is that:
\begin{align}\label{0}
du=\sum'_{|I|=p}\sum^N_{j=1}\frac{\p u_I}{\p x_j}dx_j\wedge dx^I,
\end{align}
and for $p+1$-form $f\in L^2_{p+1}(e^{-\vp})$, we say that $f$ is the differential $du$, written $du=f$, provided
\begin{align*}
\int_{G} du\cdot\a=\int_{G} f\cdot\a
\end{align*}
for all test forms $\a={\sum'_{|J|=p+1}}\a_Jdx^J\in\mathcal{D}_{p+1}$.
By the definition of ${D}_{p}$, the operator $d: \mathcal{D}_p\rightarrow \mathcal{D}_{p+1}$ is well defined.
We now extend the definition of the operator $d$ by allowing it to act on any $u\in L^2_p(e^{-\vp})$ such that $du$, in the sense of distributions, lies in $L^2_{p+1}(e^{-\vp})$. Then we obtain a closed, densely defined operator
$$T: L^2_p(e^{-\vp})\rightarrow L^2_{p+1}(e^{-\vp}),$$
where the domain of $T$ is
$$Dom(T)=\{u\in L^2_p(e^{-\vp})\mid du\in L^2_{p+1}(e^{-\vp})\}.$$

Now we consider the Hilbert space adjoint of $T$:
$$T^*: L^2_{p+1}(e^{-\vp})\rightarrow L^2_p(e^{-\vp}).$$
Let $Dom(T^*)$ be the domain of $T^*$ and $\a\in L^2_{p+1}(e^{-\vp})$. By functional analysis, we say that $\a\in Dom(T^*)$ if there exists a constant $c=c(\a)>0$ such that $$|\langle Tu,\a\rangle_{L^2_{p+1}(e^{-\vp})}|\leq c\|u\|_{L^2_p(e^{-\vp})}$$
for all $u\in Dom(T)$. This definition is equivalent to that  $\a\in Dom(T^*)$ if and only if there exists $v\in L^2_p(e^{-\vp})$ such that
$$\langle u,v\rangle_{L^2_p(e^{-\vp})}=\langle Tu,\a\rangle_{L^2_{p+1}(e^{-\vp})}$$
for all $u\in Dom(T)$. Note that $v$ is unique. We set $v=T^*\a$. Then
$T^*: Dom(T^*)\rightarrow L^2_p(e^{-\vp})$ is a linear operator and satisfies
\begin{align}\label{1}
\langle u,T^*\a\rangle_{L^2_p(e^{-\vp})}=\langle Tu,\a\rangle_{L^2_{p+1}(e^{-\vp})}
\end{align}
for all $u\in Dom(T)$, $\a\in Dom(T^*)$. It is well-known that $T^*$ is again a closed, densely defined operator.

In order to compute $T^*$, we first
compute $T^*_{formal}$, the
formal adjoint of $T$, which is defined using only test forms, i.e., we demand
\begin{align}\label{10000}
\langle Tu,\a\rangle_{L^2_{p+1}(e^{-\vp})}=\langle u,T^*_{formal}\a\rangle_{L^2_p(e^{-\vp})}
\end{align}
for $u\in Dom(T)$ and $\a\in \mathcal{D}_{p+1}$.
For $\a=\sum'_{|J|=p+1}\a_Jdx^J$, if $J_1$ is a permutation of $J$, we write $\a_{J_1}=\epsilon^{J_1}_J\a_J$, where $\epsilon^{J_1}_J$ is the signature of the permutation (for example, the signature is $-1$ if only two indices are interchanged). In particular, a term $\a_{jI}=0$ if $j\in I$, where $I=(i_1,\cdots,i_p)$ and $jI=(j,i_1,\cdots,i_p)$. Then by (\ref{0}) and
integration by parts, the left side of (\ref{10000}) is given by
\begin{align*}
\langle Tu,\a\rangle_{L^2_{p+1}(e^{-\vp})}&=\int_{G}
du\cdot\a e^{-\vp}\\
&=\int_{G}\sum'_{|I|=p}\sum^N_{j=1}\frac{\p u_I}{\p x_j}\a_{jI}e^{-\vp}\\
&=-\int_{G}\sum'_{|I|=p}\sum^N_{j=1}u_I\frac{\p(\a_{jI}e^{-\vp})}{\p x_j}\\
&=\int_{G}\lt(\sum'_{|I|=p}u_I\lt(-e^{\vp}
\sum^N_{j=1}\frac{\p(\a_{jI}e^{-\vp})}{\p x_j}\rt)\rt)e^{-\vp}\\
&=\int_{G}\lt({\sum_{|I|=p}}'u_IA_I\rt)e^{-\vp}
,
\end{align*}
where
\begin{align}\label{5}
A_I=-e^{\vp}
\sum^N_{j=1}\frac{\p(\a_{jI}e^{-\vp})}{\p x_j}.
\end{align}
For example, if $p=1$, then
$$A_I=A_i=
-e^{\vp}\lt(\sum_{1\leq j<i}\frac{\p(\a_{ji}e^{-\vp})}{\p x_j}-\sum_{i<j\leq N}\frac{\p(\a_{ij}e^{-\vp})}{\p x_j}\rt).$$
Note that $\a\in \mathcal{D}_{p+1}$ in (\ref{5}). Then $A_I$ is a smooth function with compact support in $G$, So ${\sum'_{|I|=p}}A_Idx^I\in\mathcal{D}_p\subset L^2_p(e^{-\vp})$.
Thus, the formal adjoint
\begin{align}\label{6}
T^*_{formal}\a=\sum'_{|I|=p}A_Idx^I,
\end{align}
where $A_I$ is as (\ref{5}). This implies that $\mathcal{D}_{p+1}\subset Dom(T^*)$.

In the sense of distributions, the formal adjoint $T^*_{formal}\a$ is actually well-defined for all $\a\in L^2_{p+1}(e^{-\vp})$  as $\vp$ is smooth.
We claim that
\begin{align}\label{3}
T^*\a=T^*_{formal}\a\ \ \mbox{for all}\ \ \a\in Dom(T^*).
\end{align}
Indeed, if $\a\in Dom(T^*)$, then by (\ref{1}) and
$\mathcal{D}_p\in Dom(T)$, we have for all $u\in \mathcal{D}_p$,
$$\int_{G}
u\cdot T^*\a e^{-\vp}=\langle u,T^*\a\rangle_{L^2_p(e^{-\vp})}=\langle Tu,\a\rangle_{L^2_{p+1}(e^{-\vp})}=\int_{G}
u\cdot T^*_{formal}\a e^{-\vp}.$$

\section{Proof of Poincar\'e Lemma}
We first prove some lemmas. Let $0\leq p\leq N-1$.

\bl\label{lemmaifif}
For each $p+1$-form $f$ in $L_{p+1}^2(e^{-\vp})$, there exists a solution $p$-form $u$ in $L_{p}^2(e^{-\vp})$ solving the equation
$$du=f$$ in $G$, in the sense of distributions
with the norm estimate
$$\|u\|^2_{L_{p}^2(e^{-\vp})}\leq c$$
if and only if
$$|\langle f,\a\rangle_{L_{p+1}^2(e^{-\vp})}|^2\leq c\lt\|T^*\a\rt\|^2_{L_{p}^2(e^{-\vp})}\ \ \mbox{for all}\ \ \a\in \mathcal{D}_{p+1}.$$
\el

\bp
(Necessity) For all $\a\in \mathcal{D}_{p+1}$, from the definition of $T^*$ and the Cauchy-Schwarz inequality, we have
$$
|\langle f,\a\rangle_{L_{p+1}^2(e^{-\vp})}|^2=|\langle Tu,\a\rangle_{L_{p+1}^2(e^{-\vp})}|^2=\lt|\lt\langle u,T^*\a\rt\rangle_{L_{p}^2(e^{-\vp})}\rt|^2
\leq\|u\|^2_{L_{p}^2(e^{-\vp})}
\lt\|T^*\a\rt\|^2_{L_{p}^2(e^{-\vp})}.
$$
Note that $\|u\|^2_{L_{p}^2(e^{-\vp})}\leq c$. Then $|\langle f,\a\rangle_{L_{p+1}^2(e^{-\vp})}|^2\leq c\lt\|T^*\a\rt\|^2_{L_{p}^2(e^{-\vp})}$.

(Sufficiency) Consider the subspace
$$E=\lt\{T^*\a\mid\a\in \mathcal{D}_{p+1}\rt\}\subset L_{p}^2(e^{-\vp}).$$
Define a linear functional $L_f: E\rightarrow\mb{R}$ by
$$L_f\lt(T^*\a\rt)=\langle f,\a\rangle_{L_{p+1}^2(e^{-\vp})}
.$$
Since
$$\lt|L_f\lt(T^*\a\rt)\rt|=\lt|\langle f,\a\rangle_{L_{p+1}^2(e^{-\vp})}\rt|
\leq\sqrt{c}\lt\|T^*\a\rt\|_{L_{p}^2(e^{-\vp})},$$
then $L_f$ is a bounded functional on $E$. So by Hahn-Banach's extension theorem, $L_f$ can be extended to a linear functional $\widetilde{L}_f$ on $L_{p}^2(e^{-\vp})$
such that
\begin{equation}\label{31}
\lt|\widetilde{L}_f(g)\rt|\leq\sqrt{c}\lt\|g\rt\|
_{L_{p}^2(e^{-\vp})}\ \ \mbox{for all}\ \ g\in L_{p}^2(e^{-\vp}).
\end{equation}
Using the Riesz representation theorem for $\widetilde{L}_f$, there exists a unique $u_0\in L_{p}^2(e^{-\vp})$ such that
\begin{equation}\label{32}
\widetilde{L}_f(g)=\langle u_0,g\rangle_{L_{p}^2(e^{-\vp})}\ \ \mbox{for all}\ \ g\in L_{p}^2(e^{-\vp}).
\end{equation}

Now we prove $du_0=f$. For all $\a\in \mathcal{D}_{p+1}$, apply $g=T^*\a$ in (\ref{32}). Then
$$\widetilde{L}_f\lt(T^*\a\rt)=\lt\langle u_0,T^*\a\rt\rangle_{L_{p}^2(e^{-\vp})}=\lt\langle Tu_0,\a\rt\rangle_{L_{p+1}^2(e^{-\vp})}.$$
Note that
$$\widetilde{L}_f\lt(T^*\a\rt)=L_f\lt(T^*\a\rt)=\langle f,\a\rangle_{L_{p+1}^2(e^{-\vp})}.$$
Therefore,
$$\lt\langle Tu_0,\a\rt\rangle_{L_{p+1}^2(e^{-\vp})}=\langle f,\a\rangle_{L_{p+1}^2(e^{-\vp})}\ \ \mbox{for all}\ \ \a\in \mathcal{D}_{p+1}.$$
Thus, $Tu_0=f$, i.e., $du_0=f$.

Next we give a bound for the norm of $u_0$. Let $g=u_0$ in (\ref{31}) and (\ref{32}). Then we have
$$\|u_0\|^2_{L_{p}^2(e^{-\vp})}=\lt|\langle u_0,u_0\rangle_{L_{p}^2(e^{-\vp})}\rt|=\lt|\widetilde{L}_f(u_0)\rt|
\leq\sqrt{c}\lt\|u_0\rt\|_{L_{p}^2(e^{-\vp})}.$$
Therefore, $\|u_0\|_{L_{p}^2(e^{-\vp})}^2\leq c$.

Let $u=u_0$. Then the lemma is proved.
\ep

From Lemma \ref{lemmaifif}, for proving Poincar\'e Lemma, we only need to prove that for all $\a\in \mathcal{D}_{p+1}$,
$$|\langle f,\a\rangle_{L_{p+1}^2(e^{-\vp})}|^2\leq \frac{1}{c(p+1)}\lt\|T^*\a\rt\|^2_{L_{p}^2(e^{-\vp})},$$
where $\vp, p, f, c,$ are as Poincar\'e Lemma. For this purpose, we prove the following lemmas.

\bl\label{bianjie}
Let $G=\{x\in\mathbb{R}^N\mid \rho(x)<0\}$ be a convex bounded domain, where $\rho$ is a smooth defining function. Suppose $\a=\sum'_{|J|=p+1}\a_Jdx^J$ is a smooth $p+1$-form on $\overline{G}$, and that $\a\in Dom(T^*)$. Then for any strictly increasing multiindex $I$ with $|I|=p$, we have
\begin{align*}
\sum^N_{j,k=1}\a_{kI}\frac{\p\a_{jI}}{\p x_k}\frac{\p\rho}{\p x_j}=-\sum^N_{j,k=1}\a_{jI}\a_{kI}\frac{\p^2\rho}{\p x_j\p x_k}.
\end{align*}
\el

\bp
Let $u$ be a smooth $p$-form on $\overline{G}$. Then
using integration by parts, we have
\begin{align*}
\langle Tu,\a\rangle_{L^2_{p+1}(e^{-\vp})}&=\int_G
du\cdot\a e^{-\vp}\\
&=\int_G\sum'_{|I|=p}\sum^N_{j=1}\frac{\p u_I}{\p x_j}\a_{jI}e^{-\vp}\\
&=-\int_G\sum'_{|I|=p}\sum^N_{j=1}u_I\frac{\p(\a_{jI}e^{-\vp})}{\p x_j}
+\int_{\p G}\sum'_{|I|=p}\sum^N_{j=1}u_I\a_{jI}e^{-\vp}\frac{\p\rho}{\p x_j}\frac{dS}{|\p\rho|}.
\end{align*}
Note that $$\langle Tu,\a\rangle_{L^2_{p+1}(e^{-\vp})}=\langle u,T^*\a\rangle_{L^2_p(e^{-\vp})}$$
and
$$-\int_G\sum'_{|I|=p}\sum^N_{j=1}u_I\frac{\p(\a_{jI}e^{-\vp})}{\p x_j}=\langle u,T^*\a\rangle_{L^2_p(e^{-\vp})}.$$
Then
$$\int_{\p G}\sum'_{|I|=p}\sum^N_{j=1}u_I\a_{jI}e^{-\vp}\frac{\p\rho}{\p x_j}\frac{dS}{|\p\rho|}=0$$ for any $u$.
Thus, for any strictly increasing multiindex $I$ with $|I|=p$, we have
$$\sum^N_{j=1}\a_{jI}\frac{\p\rho}{\p x_j}=0\ \ \mbox{on}\ \ \p G.$$
For the multiindex $I$ above, let $F_I=\sum^N_{j=1}\a_{jI}\frac{\p\rho}{\p x_j}$ and $L_I=\sum^N_{k=1}\a_{kI}\frac{\p}{\p x_k}$, a tangential differential operator. Then on $\p G$, we have
\begin{align*}
0=L_I(F_I)&=\sum^N_{k=1}\a_{kI}\frac{\p}{\p x_k}\lt(\sum^N_{j=1}\a_{jI}\frac{\p\rho}{\p x_j}\rt)\\
&=\sum^N_{j,k=1}\a_{kI}\frac{\p}{\p x_k}\lt(\a_{jI}\frac{\p\rho}{\p x_j}\rt)\\
&=\sum^N_{j,k=1}\a_{kI}\frac{\p\a_{jI}}{\p x_k}\frac{\p\rho}{\p x_j}+\sum^N_{j,k=1}\a_{jI}\a_{kI}\frac{\p^2\rho}{\p x_j\p x_k}
.
\end{align*}
Therefore, the lemma is proved.
\ep

\bl\label{da}
Let $\a=\sum'_{|J|=p+1}\a_Jdx^J$. Then
\begin{align}\label{dam}
|d\a|^2=\sum'_{|J|=p+1}\sum^N_{j=1}\lt|\frac{\p\a_J}{\p x_j}\rt|^2-\sum'_{|I|=p}\sum^N_{j,k=1}\frac{\p\a_{kI}}{\p x_j}\frac{\p\a_{jI}}{\p x_k}.
\end{align}
\el

\bp
Note that $d\a=0$ when $p+1=N$. When $p+1<N$,
\begin{align*}
d\a&=\sum'_{|J|=p+1}\sum^N_{j=1}\frac{\p\a_J}{\p x_j}dx_j\wedge dx^J\\
&=\sum'_{|J|=p+1}\sum_{j\notin J}\frac{\p \a_J}{\p x_j}\epsilon^{jJ}_{(jJ)'}dx^{(jJ)'}\\
&=\sum'_{|M|=p+2}\lt(\sum_{j\in M}\frac{\p \a_{M^j}}{\p x_j}\epsilon^{jM^j}_{M}\rt)dx^M,
\end{align*}
where$(jJ)'$ is the permutation of $jJ$ such that $(jJ)'$ is a strictly increasing multiindex, $\epsilon^{jJ}_{(jJ)'}$ is the signature of the permutation  and $M^j$ is the increasing multiindex with $j$ removed from $M$.
Then we prove the lemma by two cases.

Case 1: $p+1=N.$
Recall that $\a_{jI}=0$ if $j\in I$. Then for the second term on the right side of (\ref{dam}), we have
$$
\sum'_{|I|=p}\sum^N_{j,k=1}\frac{\p\a_{kI}}{\p x_j}\frac{\p\a_{jI}}{\p x_k}
=\sum'_{|I|=N-1}\sum_{j\notin I}\lt|\frac{\p\a_{jI}}{\p x_j}\rt|^2=\sum'_{|I|=N-1}\sum_{j\notin I}\lt|\frac{\p\a_{(jI)'}}{\p x_j}\rt|^2=\sum'_{|J|=N}\sum^N_{j=1}\lt|\frac{\p\a_{J}}{\p x_j}\rt|^2,
$$
which is the same as the first term on the right side of (\ref{dam}). Then (\ref{dam}) is proved for Case 1.

Case 2: $p+1<N.$ We have
\begin{align*}
|d\a|^2&=\sum'_{|M|=p+2}\lt(\sum_{j\in M}\frac{\p\a_{M^j}}{\p x_j}\epsilon^{jM^j}_M\rt)^2\\
&=\sum'_{|M|=p+2}\sum_{j\in M}\lt|\frac{\p\a_{M^j}}{\p x_j}\rt|^2
+\sum'_{|M|=p+2}\sum_{\substack{j,k\in M\\j\neq k}}\frac{\p\a_{M^j}}{\p x_j}\frac{\p\a_{M^k}}{\p x_k}\epsilon^{jM^j}_M\epsilon^{kM^k}_M\\
&=\sum'_{|J|=p+1}\sum_{j\notin J}\lt|\frac{\p\a_J}{\p x_j}\rt|^2
+\sum'_{|I|=p}\sum_{\substack{j,k\notin I\\j\neq k}}\frac{\p\a_{(kI)'}}{\p x_j}\frac{\p\a_{(jI)'}}{\p x_k}\epsilon^{j(kI)'}_{(j(kI)')'}\epsilon^{k(jI)'}_{(k(jI)')'}\\
&=\sum'_{|J|=p+1}\sum_{j\notin J}\lt|\frac{\p\a_J}{\p x_j}\rt|^2
+\sum'_{|I|=p}\sum_{\substack{j,k\notin I\\j\neq k}}\frac{\p\a_{kI}}{\p x_j}\frac{\p\a_{jI}}{\p x_k}\epsilon^{kI}_{(kI)'}\epsilon^{jI}_{(jI)'}
\epsilon^{j(kI)'}_{(j(kI)')'}\epsilon^{k(jI)'}_{(k(jI)')'}\\
&=\sum'_{|J|=p+1}\sum_{j\notin J}\lt|\frac{\p\a_J}{\p x_j}\rt|^2
-\sum'_{|I|=p}\sum_{\substack{j,k\notin I\\j\neq k}}\frac{\p\a_{kI}}{\p x_j}\frac{\p\a_{jI}}{\p x_k}.
\end{align*}
Note that
\begin{align*}
\sum'_{|J|=p+1}\sum_{j\in J}\lt|\frac{\p\a_J}{\p x_j}\rt|^2
=\sum'_{|I|=p}\sum_{j\notin I}\lt|\frac{\p\a_{(jI)'}}{\p x_j}\rt|^2
=\sum'_{|I|=p}\sum_{j\notin I}\lt|\frac{\p\a_{jI}}{\p x_j}\rt|^2.
\end{align*}
Then
\begin{align*}
|d\a|^2&=\lt(\sum'_{|J|=p+1}\sum_{j\notin J}\lt|\frac{\p\a_J}{\p x_j}\rt|^2
+\sum'_{|J|=p+1}\sum_{j\in J}\lt|\frac{\p\a_J}{\p x_j}\rt|^2\rt)-\lt(\sum'_{|I|=p}\sum_{j\notin I}\lt|\frac{\p\a_{jI}}{\p x_j}\rt|^2+\sum'_{|I|=p}\sum_{\substack{j,k\notin I\\j\neq k}}\frac{\p\a_{kI}}{\p x_j}\frac{\p\a_{jI}}{\p x_k}\rt)\\
&=\sum'_{|J|=p+1}\sum^N_{j=1}\lt|\frac{\p\a_J}{\p x_j}\rt|^2
-\sum'_{|I|=p}\sum_{j,k\notin I}\frac{\p\a_{kI}}{\p x_j}\frac{\p\a_{jI}}{\p x_k}\\
&=\sum'_{|J|=p+1}\sum^N_{j=1}\lt|\frac{\p\a_J}{\p x_j}\rt|^2
-\sum'_{|I|=p}\sum^N_{j,k=1}\frac{\p\a_{kI}}{\p x_j}\frac{\p\a_{jI}}{\p x_k}.
\end{align*}
Then (\ref{dam}) is proved for Case 2.
\ep

\bl\label{zy}
Let $G=\{x\in\mathbb{R}^N\mid \rho(x)<0\}$ be a convex bounded domain, where $\rho$ is a smooth defining function. Suppose $\a=\sum'_{|J|=p+1}\a_Jdx^J$ is a smooth $p+1$-form on $\overline{G}$, and that $\a\in Dom(T^*)$. Then
\begin{align}\label{10003}
\|T^*\a\|^2_{L^2_p(e^{-\vp})}+&\|d\a\|^2_{L^2_{p+2}(e^{-\vp})}
=\int_G\sum'_{|I|=p}\sum^N_{j,k=1}
\frac{\p^2\vp}{\p x_j\p x_k}\a_{jI}\a_{kI}e^{-\vp}\nonumber\\
&+\int_G
\sum'_{|J|=p+1}\sum^N_{j=1}\lt|\frac{\p\a_J}{\p x_j}\rt|^2 e^{-\vp}+\int_{\p G}\sum'_{|I|=p}\sum^N_{j,k=1}\a_{jI}\a_{kI}
\frac{\p^2\rho}{\p x_j\p x_k}e^{-\vp}\frac{dS}{|\p\rho|}.
\end{align}
In particular, if $G$ is a strictly convex bounded domain, and for $\vp$, there exists a constant $c>0$ such that
$$\sum^N_{j,k=1}
\frac{\p^2\vp}{\p x_j\p x_k}\xi_{j}\xi_{k}\geq c|\xi|^2\ \ \mbox{for all}\ \ \xi=(\xi_1,\cdots,\xi_N)\in\mathbb{R}^N.$$
Then
\begin{align}\label{10004}
\|T^*\a\|^2_{L^2_p(e^{-\vp})}+\|d\a\|^2_{L^2_{p+2}(e^{-\vp})}
\geq c(p+1)\|\a\|^2_{L^2_{p+1}(e^{-\vp})}.
\end{align}
\el

\bp
Consider the expression
\begin{align}\label{34}
Q=\|T^*\a\|^2_{L^2_p(e^{-\vp})}=\langle T^*\a,T^*\a\rangle_{L^2_p(e^{-\vp})}=\langle TT^*\a,\a\rangle_{L^2_p(e^{-\vp})}.
\end{align}
By (\ref{5})-(\ref{3}), we have
$$TT^*\a=d(T^*\a)=d
\lt(\sum'_{|I|=p}A_Idx^I\rt)
=\sum'_{|I|=p}\sum^N_{k=1}\frac{\p A_I}{\p x_k}dx_k\wedge dx^I,
$$
where $A_I$ is as (\ref{5}).
Let $$\delta_j=e^\vp\frac{\p}{\p x_j}e^{-\vp}
=\frac{\p}{\p x_j}-\frac{\p\vp}{\p x_j}.$$
Then
$$A_I=-\sum^N_{j=1}\delta_j\a_{jI}.$$
Observe that
$$\delta_j\frac{\p}{\p x_k}=\frac{\p^2}
{\p x_j\p x_k}-\frac{\p\vp}{\p x_j}\frac{\p}
{\p x_k}.$$
We have
$$\frac{\p}{\p x_k}\delta_j=\frac{\p^2}
{\p x_j\p x_k}-\frac{\p\vp}{\p x_j}\frac{\p}
{\p x_k}-\frac{\p^2\vp}
{\p x_j\p x_k}=\delta_j\frac{\p}{\p x_k}-\vp_{jk}.$$
Here $\vp_{jk}=\frac{\p^2\vp}
{\p x_j\p x_k}.$
So for $1\leq k\leq N$, we have
\begin{align*}
\frac{\p A_I}{\p x_k}=-\sum^N_{j=1}\lt(\frac{\p}{\p x_k}\delta_j\rt)\a_{jI}
=\sum^N_{j=1}\lt(\vp_{jk}\a_{jI}-\delta_j\frac{\p\a_{jI}}{\p x_k}\rt).
\end{align*}
Then by (\ref{34}), we have
\begin{align}
Q&=\int_G TT^*\a\cdot\a e^{-\vp}\nonumber\\
&=\int_G\sum'_{|I|=p}\sum^N_{j=1}\frac{\p A_I}{\p x_k}\a_{kI}e^{-\vp}\nonumber\\
&=\int_G\sum'_{|I|=p}\sum^N_{j,k=1}
\vp_{jk}\a_{jI}\a_{kI}e^{-\vp}
+\int_G\sum'_{|I|=p}\sum^N_{j,k=1}(-1)
\lt(\delta_j\frac{\p\a_{jI}}{\p x_k}\rt)\a_{kI}e^{-\vp}
\nonumber\\
&=Q_1+Q_2.\label{10001}
\end{align}
For $Q_2$, by Lemma \ref{bianjie} and \ref{da},we obtain that
\begin{align}
Q_2&=\int_G\sum'_{|I|=p}\sum^N_{j,k=1}(-1)\frac{\p}{\p x_j}
\lt(e^{-\vp}\frac{\p\a_{jI}}{\p x_k}\rt)\a_{kI}\nonumber\\
&=\int_G\sum'_{|I|=p}\sum^N_{j,k=1}\frac{\p\a_{kI}}{\p x_j}
\frac{\p\a_{jI}}{\p x_k}e^{-\vp}
-\int_{\p G}\sum'_{|I|=p}\sum^N_{j,k=1}\a_{kI}\frac{\p\a_{jI}}{\p x_k}e^{-\vp}\frac{\p\rho}{\p x_j}\frac{dS}{|\p\rho|}\nonumber\\
&=\int_G\sum'_{|I|=p}\sum^N_{j,k=1}\frac{\p\a_{kI}}{\p x_j}
\frac{\p\a_{jI}}{\p x_k}e^{-\vp}
+\int_{\p G}\sum'_{|I|=p}\sum^N_{j,k=1}\a_{jI}\a_{kI}
\frac{\p^2\rho}{\p x_j\p x_k}e^{-\vp}\frac{dS}{|\p\rho|}\nonumber\\
&=\int_G\sum'_{|J|=p+1}\sum^N_{j=1}\lt|\frac{\p\a_J}{\p x_j}\rt|^2 e^{-\vp}-\|d\a\|^2_{L^2_{p+2}(e^{-\vp})}+\int_{\p G}\sum'_{|I|=p}\sum^N_{j,k=1}\a_{jI}\a_{kI}
\frac{\p^2\rho}{\p x_j\p x_k}e^{-\vp}\frac{dS}{|\p\rho|}.\label{10002}
\end{align}
Then (\ref{10003}) is proved by (\ref{34}), (\ref{10001}) and (\ref{10002}).

Now we prove (\ref{10004}). Observe that
\begin{align*}
\sum'_{|I|=p}\sum^N_{j,k=1}
\vp_{jk}\a_{jI}\a_{kI}\geq\sum'_{|I|=p}c\sum^N_{j=1}\lt|\a_{jI}\rt|^2
=c\sum'_{|I|=p}\sum_{j\notin I}\lt|\a_{(jI)'}\rt|^2=c\sum'_{|J|=p+1}\sum_{j\in J}\lt|\a_{J}\rt|^2
=c(p+1)|\a|^2.
\end{align*}
Then for the first term on the right side of (\ref{10003}),
\begin{align*}
\int_G\sum'_{|I|=p}\sum^N_{j,k=1}
\vp_{jk}\a_{jI}\a_{kI}e^{-\vp}\geq\int_G c(p+1)|\a|^2e^{-\vp}
=c(p+1)\|\a\|^2_{L^2_{p+1}(e^{-\vp})}.
\end{align*}
Since $G$ is a strictly convex bounded domain, we have
$$\sum^N_{j,k=1}
\frac{\p^2\rho}{\p x_j\p x_k}\xi_{j}\xi_{k}\geq \widetilde{c}|\xi|^2\ \ \mbox{for all}\ \  \xi=(\xi_1,\cdots,\xi_N)\in\mathbb{R}^N,$$
where $\widetilde{c}$ is a positive valued function in $G$.
Then the last term on the right side of (\ref{10003}) is nonnegative.
Note that the second terms on the right side of (\ref{10003}) is always nonnegative.
Thus, (\ref{10004}) is proved.
\ep

For the proof of Poincar\'e Lemma, we need the following density lemma since the elements in $Dom(T^*)\cap Dom(S)$ are not necessarily smooth forms in general, and in the lemmas above, the computation is all based on the smooth elements.

\bl\label{bijin}
Let $f\in Dom(T^*)\cap Dom(S)$. Then there exists a sequence $\{f_\nu\}$ of smooth $p+1$-forms on $\overline{G}$, such that $f_\nu\in Dom(T^*)\cap Dom(S)$, $f_\nu\rightarrow f$ in $L^2_{p+1}(e^{-\vp})$, $T^*f_\nu\rightarrow T^*f$ in $L^2_p(e^{-\vp})$ and $Sf_\nu\rightarrow Sf$ in $L^2_{p+2}(e^{-\vp})$.
\el

The proof of this lemma would be, in principle, similar to Berndtsson's \cite {B}(Proposition 1.5.3) for his proof of the H\"{o}rmander's theorem for $\overline{\p}$, which is rather technical and nontrivial. We feel that if we had included the proof, it would have made this paper rather long. For interested readers, refer for the proof of Proposition 1.5.3 in \cite {B}.

Now we give the proof of Poincar\'e Lemma.

\bp
Let $N=\{f\mid f\in L^2_{p+1}(e^{-\vp});df=0\}$, which is a closed subspace of $L^2_{p+1}(e^{-\vp})$.
For each $\a$ in $\mathcal{D}_{p+1}$, clearly $\a\in L^2_{p+1}(e^{-\vp})$, so we can decompose $\a=\a^1+\a^2$, where $\a^1$ lies in $N$ and $\a^2$ is orthogonal to $N$. This implies that $\a^2$ is orthogonal to any form $Tu$, since $Tu\in N$. So by the definition of $Dom(T^*)$, we see that $\a^2$ lies in the domain of $T^*$ and $T^*\a^2=0$. Since $\a$ lies in the domain of $T^*$, it follows that $T^*\a=T^*\a^1$.

Note that $\a^1\in Dom(T^*)\cap Dom(S)$. Then by Lemma \ref{bijin},
there exists a sequence $\{\a_\nu\}$, which are smooth $p+1$-forms on $\overline{G}$, such that $\a_\nu\in Dom(T^*)\cap Dom(S)$, $\a_\nu\rightarrow\a^1$ in $L_{p+1}^2(e^{-\vp})$, $T^*\a_\nu\rightarrow T^*\a^1$ in $L_{p}^2(e^{-\vp})$, and $S\a_\nu\rightarrow S\a^1$ in $L_{p+2}^2(e^{-\vp})$.

For $\a_\nu$, by Lemma \ref{zy}, we have
\begin{align*}
\|T^*\a_\nu\|^2_{L^2_p(e^{-\vp})}+\|S\a_\nu\|^2_{L^2_{p+2}(e^{-\vp})}\geq
c(p+1)\|\a_\nu\|^2_{L^2_{p+1}(e^{-\vp})}.
\end{align*}
Let $\nu\rightarrow+\infty$, so
\begin{align*}
\|T^*\a^1\|^2_{L^2_p(e^{-\vp})}+\|S\a^1\|^2_{L^2_{p+2}(e^{-\vp})}\geq
c(p+1)\|\a^1\|^2_{L^2_{p+1}(e^{-\vp})},
\end{align*}
which means that
\begin{align*}
\|T^*\a^1\|^2_{L^2_p(e^{-\vp})}\geq
c(p+1)\|\a^1\|^2_{L^2_{p+1}(e^{-\vp})}
\end{align*}
since $S\a^1=0$.

By the Cauchy-Schwarz inequality, we have
\begin{align*}
\lt|\langle f,\a^1\rangle_{L_{p+1}^2(e^{-\vp})}\rt|^2
&\leq\lt\|f\rt\|^2_{L_{p+1}^2(e^{-\vp})}
\lt\|\a^1\rt\|^2_{L_{p+1}^2(e^{-\vp})}\\
&=\lt(\frac{1}{c(p+1)}\lt\|f\rt\|^2_{L_{p+1}^2(e^{-\vp})}\rt)
\lt(c(p+1)\lt\|\a^1\rt\|^2_{L_{p+1}^2(e^{-\vp})}\rt)\\
&\leq\lt(\frac{1}{c(p+1)}\lt\|f\rt\|^2_{L_{p+1}^2(e^{-\vp})}\rt)
\|T^*\a^1\|^2_{L^2_p(e^{-\vp})}.
\end{align*}
Let $\hat{c}=\frac{1}{c(p+1)}\lt\|f\rt\|^2_{L_{p+1}^2(e^{-\vp})}$. Then
\begin{align*}
\lt|\langle f,\a^1\rangle_{L_{p+1}^2(e^{-\vp})}\rt|^2
\leq \hat{c}
\|T^*\a^1\|^2_{L^2_p(e^{-\vp})}\ \ \mbox{for all}\ \ \a\in \mathcal{D}_{p+1}.
\end{align*}
Note that $f\in N$.
Thus,
$$\lt|\langle f,\a\rangle_{L_{p+1}^2(e^{-\vp})}\rt|^2=\lt|\langle f,\a^1\rangle_{L_{p+1}^2(e^{-\vp})}\rt|^2\leq \hat{c}\lt\|T^*\a^1\rt\|^2_{L_{p}^2(e^{-\vp})}
=\hat{c}\lt\|T^*\a\rt\|^2_{L_{p}^2(e^{-\vp})}.$$
By Lemma \ref{lemmaifif}, there exists a solution $u\in L_{p}^2(e^{-\vp})$ solving the equation
$$du=f$$ in $G$
with the norm estimate
$\|u\|^2_{L_{p}^2(e^{-\vp})}\leq \hat{c}$,
i.e.,
$$\int_G|u|^2e^{-\vp}\leq
\frac{1}{c(p+1)}\int_G|f|^2e^{-\vp}.$$
The theorem is proved.
\ep

\section{Preliminary for the main theorem}

Let $n\geq1$ be an integer and $\O$ be a strictly convex bounded domain in $\mathbb{C}^{n}$. Let $\vp$ be a nonnegative smooth function on $\overline{\O}$ and the weighted Hilbert space
$$L^2(\O,e^{-\vp})
=\{u:\O\rightarrow\mb{C}\mid u\in L^2_{loc}(\O); \int_{\O}|u|^2 e^{-\vp}<+\infty\}.$$
We denote
the weighted inner product for $u,v\in L^2(\O,e^{-\vp})$ by
$$\langle u,v\rangle_{L^2(\O,e^{-\vp})}=\int_{\O}u\overline{v} e^{-\vp},$$
and the weighted norm of $u\in L^2(\O,e^{-\vp})$ by
$\|u\|_{L^2(\O,e^{-\vp})}=\sqrt{\langle u,u\rangle}_{L^2(\O,e^{-\vp})}.$

In general, a $(1,1)$ form $f$ on $\O$ is a formal combination
$$f=\sum^n_{i,j=1}f_{i\overline{j}}dz_i\wedge d\overline{z}_j,$$
where $f_{i\overline{j}}:\O\rightarrow\mb{C}$ is a function for $1\leq i,j\leq n$. For $(1,1)$ forms $f$ and $g$, we denote their pointwise scalar product by
$$f\cdot \overline{g}=\sum^n_{i,j=1}f_{i\overline{j}}\overline{g}_{i\overline{j}}.$$
We also consider the weighted Hilbert space for $(1,1)$ forms
$$L^2_{(1,1)}(\O,e^{-\vp})
=\{f=\sum^n_{i,j=1}f_{i\overline{j}}dz_i\wedge d\overline{z}_j\mid f_{i\overline{j}}\in L^2_{loc}(\O); \int_{\O}|f|^2 e^{-\vp}<+\infty\},$$
where $|f|^2=f\cdot \overline{f}$. We denote
the weighted inner product for $f,g\in L^2_{(1,1)}(\O,e^{-\vp})$ by
$$\langle f,g\rangle_{L^2_{(1,1)}(\O,e^{-\vp})}=\int_{\O}
f\cdot \overline{g} e^{-\vp},$$
and the weighted norm of $f\in L^2_{(1,1)}(\O,e^{-\vp})$ by
$\|f\|_{L^2_{(1,1)}(\O,e^{-\vp})}=\sqrt{\langle f,f\rangle}_{L^2_{(1,1)}(\O,e^{-\vp})}.$

For the conversion between complex and real forms, we need the following lemmas, which can be verified by simple computations.

\bl\label{wj1}
Let $f=\sum^n_{i,j=1}(A_{i\overline{j}}+\sqrt{-1}B_{i\overline{j}})dz_i\wedge d\overline{z}_j$ be any real $(1,1)$ form (i.e. $f=\overline{f}$), where $A_{i\overline{j}}$ and $B_{i\overline{j}}$ are real-valued functions. Then $A_{i\overline{j}}=-A_{j\overline{i}}, B_{i\overline{j}}=B_{j\overline{i}}$ and
$f$ can be decomposed to a real $2$-form
$$f=2\lt(\sum_{1\leq i<j\leq n}A_{i\overline{j}}dx_i\wedge dx_j
+\sum_{1\leq i<j\leq n}A_{i\overline{j}}dy_i\wedge dy_j +\sum^n_{i,j=1}B_{i\overline{j}}
dx_i\wedge dy_j\rt).$$
\el

\bl\label{wj2}
Let $v=\sum^{2n}_{j=1}v_jdx_j$ be any real $1$-form, where $v_j$ and are real-valued functions. Let $z_j=x_{2j-1}+\sqrt{-1}x_{2j}$. Then
$v$ can be decomposed to $$v=v^{1,0}+v^{0,1},$$
where $v^{1,0}=\sum^{n}_{j=1}\lt(\frac{1}{2}v_{2j-1}
+\frac{1}{2\sqrt{-1}}v_{2j}\rt)dz_j$ is a $(1,0)$ form and
$v^{0,1}=\sum^{n}_{j=1}\lt(\frac{1}{2}v_{2j-1}
-\frac{1}{2\sqrt{-1}}v_{2j}\rt)d\overline{z}_j$ is a $(0,1)$ form.
\el

\section{Proof of the main theorem}
First we give three lemmas. They are all well-known and can be simply verified by virtue of the definition of distributions.

\bl\label{pu}
If $u\in L^2(\O,e^{-\vp})$ and $\overline{\p}u\in L^2_{{0,1}}(\O,e^{-\vp})$. Then $\p\overline{u}=\overline{\overline{\p}u}$, where $\overline{\p}u$ and $\p\overline{u}$ are in the sense of distributions.
\el

\bl\label{huhuan}
If $u\in L^2(\O,e^{-\vp})$. Then $\overline{\p}\p u=-\p\overline{\p}u$ in the sense of distributions.
\el

\noindent\textbf{Remark 5.1.} In the lemma, it is crucial that $\overline{\p}\p u$ and $\p\overline{\p}u$ are both forms. Otherwise, when $n=1$, $\overline{\p}\p u=\p\overline{\p}u$ if $\overline{\p}\p u=\p\overline{\p}u=\frac{\p^2u}{\p z\p \overline{z}}$ are as weak derivatives.

\bl\label{zuijiao}
Let $u\in L^2(\O,e^{-\vp})$. If $\overline{\p}u\in L^2_{{0,1}}(\O,e^{-\vp})$, then $\p\overline{\p}u=\p(\overline{\p}u)$ in the sense of distributions. If $\p u\in L^2_{{1,0}}(\O,e^{-\vp})$, then $\overline{\p}\p u=\overline{\p}(\p u)$ in the sense of distributions.
\el

To prove the main theorem, we also need the following simple version of the H\"{o}rmander's theorem \cite{4} (page 92, Lemma 4.4.1).\\

\noindent\textbf{H\"{o}rmander's theorem.} (A simple version for $(0,1)$ forms)
\textit{Let $\O$ be a pseudoconvex open set in $\mathbb{C}^{n}$. Let $\vp$ be a real-valued smooth function in $\O$ such that
$$\sum^n_{j,k=1}
\frac{\p^2\vp}{\p z_j\p \overline{z}_k}\omega_{j}\overline{\omega}_{k}\geq c|\omega|^2\ \ \mbox{for all}\ \  \omega=(\omega_1,\cdots,\omega_n)\in\mathbb{C}^n,$$
where $c>0$ is a constant.
For each $f\in L^2_{(0,1)}(\O,e^{-\vp})$ such that $\overline{\p}f=0$, there exists a solution $u$ in $L^2(\O,e^{-\vp})$ solving equation
$\overline{\p}u=f$ in $\O$, in the sense of distributions,
with the norm estimate $$\int_{\O} |u|^2e^{-\vp}\leq
\frac{2}{c}\int_{\O}|f|^2e^{-\vp}.$$}

In order to apply the H\"{o}rmander's theorem above, we need the following results that convert real convexity to  plurisubharmonicity, which is well-known, but for which a brief proof is provided.

\bl\label{dengshi}
Let $\vp$ be a smooth function in a domain in $\mathbb{R}^{2n}$. Let $\xi=(\xi_1,\cdots,\xi_{2n})\in\mathbb{R}^{2n}$ and $x=(x_1,\cdots,x_{2n})\in\mathbb{R}^{2n}$. For $1\leq j\leq n$, let $\omega_j=\xi_j+\sqrt{-1}\xi_{j+n}$ and $z_j=x_j+\sqrt{-1}x_{j+n}$.
Then
$$\sum^{2n}_{j,k=1}
\frac{\p^2\vp}{\p x_j\p x_k}\xi_j\xi_k=\sum^n_{j,k=1}\lt(\frac{\p^2\vp}{\p z_j\p z_k}\omega_j\omega_k
+2\frac{\p^2\vp}{\p z_j\p \overline{z}_k}\omega_j
\overline{\omega}_k+
\frac{\p^2\vp}{\p \overline{z}_j\p \overline{z}_k}\overline{\omega}_j
\overline{\omega}_k\rt).$$
\el

\bp
Consider $\vp$ in real variables. Let $\phi(t)=\vp(x+t\xi)$. Then the second derivative of $\phi$ at $0$ is
$$\phi''(0)=\sum^{2n}_{j,k=1}
\frac{\p^2\vp}{\p x_j\p x_k}\xi_j\xi_k.$$
Consider the same function $\phi$ in complex variables as $\phi(t)=\vp(z+tw)$. Then the second derivative of $\phi$ at $0$ is
$$\phi''(0)=\sum^n_{j,k=1}\lt(\frac{\p^2\vp}{\p z_j\p z_k}\omega_j\omega_k
+2\frac{\p^2\vp}{\p z_j\p \overline{z}_k}\omega_j
\overline{\omega}_k+
\frac{\p^2\vp}{\p \overline{z}_j\p \overline{z}_k}\overline{\omega}_j
\overline{\omega}_k\rt).$$
Thus the lemma is proved.
\ep

From the lemma above, we have the following lemma, which we need to use.

\bl\label{convexdengshi}
Let $\vp$ be a real-valued smooth function in a domain in $\mathbb{R}^{2n}$ such that
$$\sum^{2n}_{j,k=1}
\frac{\p^2\vp}{\p x_j\p x_k}\xi_j\xi_k\geq c|\xi|^2\ \ \mbox{for all}\ \ \xi=(\xi_1,\cdots,\xi_{2n})\in\mathbb{R}^{2n},$$
where $c>0$ is a constant.
For $1\leq j\leq n$, let $z_j=x_j+\sqrt{-1}x_{j+n}$.
Then
$$
\sum^n_{j,k=1}\frac{\p^2\vp}{\p z_j\p \overline{z}_k}\omega_j
\overline{\omega}_k\geq \frac{c}{2}|\omega|^2+\lt|\sum^n_{j,k=1}\frac{\p^2\vp}{\p z_j\p z_k}\omega_j\omega_k\rt|\ \ \mbox{for all}\ \ \omega=(\omega_1,\cdots,\omega_{n})\in\mathbb{C}^{n}.
$$
\el

\bp
By Lemma \ref{dengshi}, we have
$$\sum^{2n}_{j,k=1}
\frac{\p^2\vp}{\p x_j\p x_k}\xi_{j}\xi_{k}=2\sum^n_{j,k=1}\lt(Re\lt(\frac{\p^2\vp}{\p z_j\p z_k}\omega_j\omega_k\rt)
+\frac{\p^2\vp}{\p z_j\p \overline{z}_k}\omega_j
\overline{\omega}_k\rt),$$
where $\omega_j=\xi_j+\sqrt{-1}\xi_{j+n}$ for $1\leq j\leq n$.
Then for all $ \omega=(\omega_1,\cdots,\omega_{n})\in\mathbb{C}^{n}$,
$$
2\sum^n_{j,k=1}\lt(Re\lt(\frac{\p^2\vp}{\p z_j\p z_k}\omega_j\omega_k\rt)
+\frac{\p^2\vp}{\p z_j\p \overline{z}_k}\omega_j
\overline{\omega}_k\rt)\geq c|\omega|^2,
$$
i.e.,
$$
\sum^n_{j,k=1}\frac{\p^2\vp}{\p z_j\p \overline{z}_k}\omega_j
\overline{\omega}_k\geq \frac{c}{2}|\omega|^2-Re\lt(\sum^n_{j,k=1}\frac{\p^2\vp}{\p z_j\p z_k}\omega_j\omega_k\rt).
$$
For any fixed $\omega\in\mathbb{C}^{n}$, replace $\omega$ by $e^{\sqrt{-1}\theta}\omega$ in the above formula, where $\theta$ is a real number such that $$-Re\lt(e^{2\sqrt{-1}\theta}\sum^n_{j,k=1}\frac{\p^2\vp}{\p z_j\p z_k}\omega_j\omega_k\rt)=\lt|\sum^n_{j,k=1}\frac{\p^2\vp}{\p z_j\p z_k}\omega_j\omega_k\rt|.$$
Then we have
$$
\sum^n_{j,k=1}\frac{\p^2\vp}{\p z_j\p \overline{z}_k}\omega_j
\overline{\omega}_k\geq \frac{c}{2}|\omega|^2+\lt|\sum^n_{j,k=1}\frac{\p^2\vp}{\p z_j\p z_k}\omega_j\omega_k\rt|.
$$
Thus,
the lemma is proved.
\ep

Finally, we are ready to give the proof of the main theorem.

\bp
First we prove the theorem for the case that $f$ is a real $(1,1)$ form. Observe that $\sqrt{-1}\p\overline{\p}$ is a real operator by Lemma \ref{huhuan}.

For real $(1,1)$ form $f\in L^2_{(1,1)}(\O,e^{-\vp})$, by Lemma \ref{wj1} $f$ can be seen as a real $2$-form and $f\in L^2_{2}(\O,e^{-\vp})$. Then by Poincar\'e Lemma ($p=1$), there exists $v\in L^2_{1}(\O,e^{-\vp})$ such that
\begin{align}\label{15}
dv=f
\end{align}
with
\begin{align}\label{zygs1}
\|v\|^2_{L^2_{1}(\O,e^{-\vp})}
\leq
\frac{1}{2c}\|f\|^2_{L^2_{2}(\O,e^{-\vp})}.
\end{align}
For this $v$, by Lemma \ref{wj2} we have the decomposition in complex forms
\begin{align}\label{16}
v=v^{1,0}+v^{0,1},
\end{align}
where $v^{1,0}\in L^2_{{1,0}}(\O,e^{-\vp})$, $v^{0,1}\in L^2_{{0,1}}(\O,e^{-\vp})$, $\overline{v^{1,0}}=v^{0,1}$ and $\overline{v^{0,1}}=v^{1,0}$.
By (\ref{15}) and (\ref{16}), we have
\begin{align}\label{17}
f
=(\p+\overline{\p})(v^{1,0}+v^{0,1})
=\p v^{1,0}+\p v^{0,1}+\overline{\p}v^{1,0} +\overline{\p}v^{0,1}.
\end{align}
Note that $\p v^{1,0}$ is a $(2,0)$ form, $\overline{\p}v^{0,1}$ is a $(0,2)$ form and $f$ can be seen as a $(1,1)$ form.
So from (\ref{17}), we have $\p v^{1,0}=0$, $\overline{\p}v^{0,1}=0$ and
\begin{align}\label{18}
\p v^{0,1}+\overline{\p}v^{1,0}=f.
\end{align}
For $v^{0,1}$, by Lemma \ref{convexdengshi} and the H\"{o}rmander's theorem (replace $c$ by $\frac{c}{2}$), there exists $w\in L^2(\O,e^{-\vp})$ such that
\begin{align}\label{19}
\overline{\p}w=v^{0,1},
\end{align}
with
\begin{align}\label{zygs2}
\int_{\O} |w|^2e^{-\vp}\leq
\frac{4}{c}\int_{\O}\lt|v^{0,1}\rt|^2e^{-\vp}.
\end{align}
So for $w$, by Lemma \ref{pu} and $\overline{v^{0,1}}=v^{1,0}$, we have
\begin{align}\label{20}
\p\overline{w}=\overline{\overline{\p}w}=v^{1,0}.
\end{align}
Then by Lemma \ref{wj2}, \ref{huhuan}, \ref{zuijiao}, (\ref{zygs1}), and (\ref{18})-(\ref{20}), we obtain
\begin{align}\label{21}
\p\overline{\p}\lt(w-\overline{w}\rt)=
\p\overline{\p}w-\p\overline{\p}\overline{w}=
\p\overline{\p}w+\overline{\p}\p\overline{w}
=\p(\overline{\p}w)+\overline{\p}(\p\overline{w})=\p v^{0,1}+\overline{\p}v^{1,0}=f,
\end{align}
with
\begin{align}
\int_{\O} |w-\overline{w}|^2e^{-\vp}&\leq
4\int_{\O}|w|^2e^{-\vp}
\leq\frac{16}{c}\int_{\O}\lt|v^{0,1}\rt|^2e^{-\vp}
=\frac{4}{c}\|v\|^2_{L^2_{1}(\O,e^{-\vp})}\leq
\frac{2}{c^2}\|f\|^2_{L^2_{2}(\O,e^{-\vp})}.\label{zygs3}
\end{align}
Let $u=-\sqrt{-1}\lt(w-\overline{w}\rt)$. Note that $\|f\|^2_{L^2_{2}(\O,e^{-\vp})}=4\|f\|^2_{L^2_{(1,1)}(\O,e^{-\vp})}$.
Then $u\in L^2(\O,e^{-\vp})$ and $$\sqrt{-1}\p\overline{\p}u=f\ \ \mbox{with}\ \ \|u\|^2_{L^2(\O,e^{-\vp})}\leq\frac{8}{c^2}
\|f\|^2_{L^2_{(1,1)}(\O,e^{-\vp})}.$$
So the theorem is proved for the case that $f$ is a real $(1,1)$ form.

Now we prove the theorem for the case that $f$ is not a real $(1,1)$ form.
Write $f=f_1+\sqrt{-1}f_2$, where $f_1=\frac{1}{2}(f+\overline{f})$ and $f_2=\frac{1}{2\sqrt{-1}}(f-\overline{f})$. Then $f_1$ and $f_2$ are real $2$-forms. Apply twice the same way above and then the theorem is proved.
\ep

\end{document}